\documentclass[12pt]{article}
\usepackage{amssymb}
\usepackage{amsfonts}
\usepackage{amsmath}
\usepackage{amsbsy}
\newcommand{\be}{\begin{equation}}
\newcommand{\ee}{\end{equation}}
\newcommand{\bea}{\begin{eqnarray}}
\newcommand{\eea}{\end{eqnarray}}
\newcommand{\ba}{\begin{array}}
\newcommand{\ea}{\end{array}}

\newcommand{\bc}{\begin{center}}
\newcommand{\ec}{\end{center}}
\newcommand{\ben}{\begin{enumerate}}
\newcommand{\een}{\end{enumerate}}
\newcommand{\bfi}{\begin{figure}}
\newcommand{\efi}{\end{figure}}

\newcommand{\bq}{\begin{quote}}
\newcommand{\eq}{\end{quote}}
\newcommand{\bqu}{\begin{quotation}}
\newcommand{\equ}{\end{quotation}}
\newenvironment{emphit}{\begin{itemize}}{\end{itemize}}
\newcommand{\bemp}{\begin{emphit}}
\newcommand{\eemp}{\end{emphit}}

\newcommand{\bt}{\begin{tabular}}
\newcommand{\et}{\end{tabular}}

\newtheorem{myth}{Theorem}[section]

\newtheorem{mydef}{Definition}[section]

\begin{document}
\date{}
\title{A Double Cryptography Using The Smarandache Keedwell Cross Inverse Quasigroup\footnote{2000
Mathematics Subject Classification. Primary 20NO5 ; Secondary 08A05}
\thanks{{\bf Keywords and Phrases :} Smarandache holomorph of loops, Smarandache cross inverse property
quasigroups(CIPQs), Smarandache automorphism group, cryptography}}
\author{T\`em\'it\d{\'o}p\d{\'e} Gb\d{\'o}l\'ah\`an Ja\'iy\'e\d ol\'a\thanks{All correspondence to be addressed to this author.} \\
Department of Mathematics,\\
Obafemi Awolowo University, Ile Ife, Nigeria.\\
jaiyeolatemitope@yahoo.com, tjayeola@oauife.edu.ng}
 \maketitle
\begin{abstract}
The present study further strengthens the use of the Keedwell CIPQ
against attack on a system by the use of the Smarandache Keedwell
CIPQ for cryptography in a similar spirit in which the cross inverse
property has been used by Keedwell. This is done as follows. By
constructing two S-isotopic S-quasigroups(loops) $U$ and $V$ such
that their Smarandache automorphism groups are not trivial, it is
shown that $U$ is a SCIPQ(SCIPL) if and only if $V$ is a
SCIPQ(SCIPL). Explanations and procedures are given on how these
SCIPQs can be used to double encrypt information.
\end{abstract}

\section{Introduction}
\subsection{Quasigroups And Loops}
Let $L$ be a non-empty set. Define a binary operation ($\cdot $) on
$L$ : If $x\cdot y\in L$ for all $x, y\in L$, $(L, \cdot )$ is
called a groupoid. If the system of equations ;
\begin{displaymath}
a\cdot x=b\qquad\textrm{and}\qquad y\cdot a=b
\end{displaymath}
have unique solutions for $x$ and $y$ respectively, then $(L, \cdot
)$ is called a quasigroup. For each $x\in L$, the elements $x^\rho
=xJ_\rho ,x^\lambda =xJ_\lambda\in L$ such that $xx^\rho=e^\rho$ and
$x^\lambda x=e^\lambda$ are called the right, left inverses of $x$
respectively. Now, if there exists a unique element $e\in L$ called
the identity element such that for all $x\in L$, $x\cdot e=e\cdot
x=x$, $(L, \cdot )$ is called a loop.  To every loop $(L,\cdot )$
with automorphism group $AUM(L,\cdot )$, there corresponds another
loop. Let the set $H=(L,\cdot )\times AUM(L,\cdot )$. If we define
'$\circ$' on $H$ such that $(\alpha, x)\circ (\beta,
y)=(\alpha\beta, x\beta\cdot y)$ for all $(\alpha, x),(\beta, y)\in
H$, then $H(L,\cdot )=(H,\circ)$ is a loop as shown in Bruck
\cite{phd82} and is called the Holomorph of $(L,\cdot )$.
\paragraph{}
A loop is a weak inverse property loop(WIPL) if and only if it obeys
the identity
\begin{equation*}\label{eq:8}
x(yx)^\rho=y^\rho\qquad\textrm{or}\qquad(xy)^\lambda x=y^\lambda.
\end{equation*}
A loop(quasigroup) is a cross inverse property
loop(quasigroup)[CIPL(CIPQ)] if and only if it obeys the identity
\begin{equation*}\label{eq:8.1}
xy\cdot x^\rho =y\qquad\textrm{or}\qquad x\cdot yx^\rho
=y\qquad\textrm{or}\qquad x^\lambda\cdot
(yx)=y\qquad\textrm{or}\qquad x^\lambda y\cdot x=y.
\end{equation*}
A loop(quasigroup) is an automorphic inverse property
loop(quasigroup)[AIPL(AIPQ)] if and only if it obeys the identity
\begin{equation*}
(xy)^\rho=x^\rho y^\rho~or~(xy)^\lambda =x^\lambda y^\lambda.
\end{equation*}
The set $SYM(G, \cdot )=SYM(G)$ of all bijections in a groupoid
$(G,\cdot )$ forms a group called the permutation(symmetric) group
of the groupoid $(G,\cdot )$. Consider $(G,\cdot )$ and $(H,\circ )$
been two distinct groupoids(quasigroups, loops). Let $A,B$ and $C$
be three distinct non-equal bijective mappings, that maps $G$ onto
$H$. The triple $\alpha =(A,B,C)$ is called an isotopism of
$(G,\cdot )$ onto $(H,\circ )$ if and only if
\begin{displaymath}
xA\circ yB=(x\cdot y)C~\forall~x,y\in G.
\end{displaymath}
If $(G,\cdot )=(H,\circ )$, then the triple $\alpha =(A,B,C)$ of
bijections on $(G,\cdot )$ is called an autotopism of the
groupoid(quasigroup, loop) $(G,\cdot )$. Such triples form a group
$AUT(G,\cdot )$ called the autotopism group of $(G,\cdot )$.
Furthermore, if $A=B=C$, then $A$ is called an automorphism of the
groupoid(quasigroup, loop) $(G,\cdot )$. Such bijections form a
group $AUM(G,\cdot )$ called the automorphism group of $(G,\cdot )$.

As observed by Osborn \cite{phd89}, a loop is a WIPL and an AIPL if
and only if it is a CIPL. The past efforts of Artzy \cite{phd140,
phd193, phd158, phd30}, Belousov and Tzurkan \cite{phd192} and
recent studies of Keedwell \cite{phd176}, Keedwell and Shcherbacov
\cite{phd175, phd177, phd178} are of great significance in the study
of WIPLs, AIPLs, CIPQs and CIPLs, their generalizations(i.e
m-inverse loops and quasigroups, (r,s,t)-inverse quasigroups) and
applications to cryptography.

Interestingly, Huthnance \cite{phd44} showed that if $(L,\cdot )$ is
a loop with holomorph $(H,\circ)$, $(L,\cdot )$ is a WIPL if and
only if $(H,\circ)$ is a WIPL. But the holomorphic structure of AIPL
and a CIPL has just been revealed by Ja\'iy\'e\d ol\'a \cite{sma15}.

In the quest for the application of CIPQs with long inverse cycles
to cryptography, Keedwell \cite{phd176} constructed the following
CIPQ which we shall specifically call Keedwell CIPQ.
\begin{myth}(Keedwell CIPQ)
Let $(G,\cdot )$ be an abelian group of order $n$ such that $n+1$ is
composite. Define a binary operation '$\circ$' on the elements of
$G$ by the relation $a\circ b=a^rb^s$, where $rs=n+1$. Then
$(G,\circ )$ is a CIPQ and the right crossed inverse of the element
$a$ is $a^u$, where $u=(-r)^3$
\end{myth}
The author also gave examples and detailed explanation and
procedures of the use of this CIPQ for cryptography. Cross inverse
property quasigroups have been found appropriate for cryptography
because of the fact that the left and right inverses $x^\lambda$ and
$x^\rho$ of an element $x$ do not coincide unlike in left and right
inverse property loops, hence this gave rise to what is called
'cycle of inverses' or 'inverse cycles' or simply 'cycles' i.e
finite sequence of elements $x_1,x_2,\cdots ,x_n$ such that
$x_k^\rho =x_{k+1}~\bmod{n}$. The number $n$ is called the length of
the cycle. The origin of the idea of cycles can be traced back to
Artzy \cite{phd140,phd193} where he also found there existence in
WIPLs apart form CIPLs. In his two papers, he proved some results on
possibilities for the values of $n$ and for the number $m$ of cycles
of length $n$ for WIPLs and especially CIPLs. We call these "Cycle
Theorems" for now.

In application, it is assumed that the message to be transmitted can
be represented as single element $x$ of a quasigroup $(L,\cdot )$
and that this is enciphered by multiplying by another element $y$ of
$L$ so that the encoded message is $yx$. At the receiving end, the
message is deciphered by multiplying by the right inverse $y^\rho$
of $y$. If a left(right) inverse quasigroup is used and the
left(right) inverse of $x$ is $x^\lambda$ ($x^\rho$), then the
left(right) inverse of $x^\lambda$ ($x^\rho$) is necessarily $x$.
But if a CIPQ is used, this is not necessary the situation. This
fact makes an attack on the system more difficult in the case of
CIPQs.

\subsection{Smarandache Quasigroups And Loops}
The study of Smarandache loops was initiated by W. B. Vasantha
Kandasamy in 2002. In her book \cite{phd75}, she defined a
Smarandache loop(S-loop) as a loop with at least a subloop which
forms a subgroup under the binary operation of the loop. In
\cite{sma2}, the present author defined a Smarandache
quasigroup(S-quasigroup) to be a quasigroup with at least a
non-trivial associative subquasigroup called a Smarandache
subsemigroup (S-subsemigroup). Examples of Smarandache quasigroups
are given in Muktibodh \cite{muk2}. In her book, she introduced over
75 Smarandache concepts on loops. In her first paper \cite{phd83},
on the study of Smarandache notions in algebraic structures, she
introduced Smarandache : left(right) alternative loops, Bol loops,
Moufang loops, and Bruck loops. But in \cite{sma1}, the present
author introduced Smarandache : inverse property loops(IPL) and weak
inverse property loops(WIPL).

A quasigroup(loop) is called a Smarandache "certain"
quasigroup(loop) if it has at least a non-trivial
subquasigroup(subloop) with the "certain" property and the latter is
referred to as the Smarandache "certain" subquasigroup(subloop). For
example, a loop is called a Smarandache CIPL(SCIPL) if it has at
least a non-trivial subloop that is a CIPL and the latter is
referred to as the Smarandache CIP-subloop. By an "initial
S-quasigroup" $L$ with an "initial S-subquasigroup" $L'$, we mean
$L$ and $L'$ are pure quasigroups, i.e. they do not obey a "certain"
property(not of any variety).

If $L$ is a S-groupoid with a S-subsemigroup $H$, then the set
$SSYM(L, \cdot )=SSYM(L)$ of all bijections $A$ in $L$ such that
$A~:~H\to H$ forms a group called the Smarandache
permutation(symmetric) group of the S-groupoid. In fact, $SSYM(L)\le
SYM(L)$.

\begin{mydef}\label{1:1}
Let $(L, \cdot )$ and $(G, \circ )$ be two distinct groupoids that
are isotopic under a triple $(U, V, W)$. Now, if $(L, \cdot )$ and
$(G, \circ )$ are S-groupoids with S-subsemigroups $L'$ and $G'$
respectively such that $A~:~L'\to G'$, where $A\in\{U,V,W\}$, then
the isotopism $(U, V, W) : (L, \cdot )\rightarrow (G, \circ )$ is
called a Smarandache isotopism(S-isotopism).

Thus, if $U=V=W$, then $U$ is called a Smarandache isomorphism,
hence we write $(L, \cdot )\succsim (G, \circ )$.

But if $(L, \cdot )=(G, \circ )$, then the autotopism $(U, V, W)$ is
called a Smarandache autotopism (S-autotopism) and they form a group
$SAUT(L,\cdot )$ which will be called the Smarandache autotopism
group of $(L, \cdot )$. Observe that $SAUT(L,\cdot )\le AUT(L,\cdot
)$. Furthermore, if $U=V=W$, then $U$ is called a Smarandache
automorphism of $(L,\cdot )$. Such Smarandache permutations form a
group $SAUM(L,\cdot )$ called the Smarandache automorphism
group(SAG) of $(L,\cdot )$.
\end{mydef}
Now, set $H_S=(L,\cdot )\times SAUM(L,\cdot )$. If we define
'$\circ$' on $H_S$ such that $(\alpha, x)\circ (\beta,
y)=(\alpha\beta, x\beta\cdot y)$ for all $(\alpha, x),(\beta, y)\in
H_S$, then $H_S(L,\cdot )=(H_S,\circ)$ is a S-quasigroup(S-loop)
with S-subgroup $(H',\circ )$ where $H'=L'\times SAUM(L)$ and thus
will be called the Smarandache Holomorph(SH) of $(L,\cdot )$.
\paragraph{}
The aim of the present study is to further strengthen the use of the
Keedwell CIPQ against attack on a system by the use of the
Smarandache Keedwell CIPQ for cryptography in a similar spirit in
which the cross inverse property has been used by Keedwell. This is
done as follows. By constructing two S-isotopic S-quasigroups(loops)
$U$ and $V$ such that their Smarandache automorphism groups are not
trivial, it is shown that $U$ is a SCIPQ(SCIPL) if and only if $V$
is a SCIPQ(SCIPL). Explanations and procedures are given on how
these SCIPQs can be used to double encrypt information.

\section{Preliminary Results}
\begin{mydef}(Smarandache Keedwell CIPQ)

Let $Q$ be an initial S-quasigroup with an initial S-subquasigroup
$P$. $Q$ is called a Smarandache Keedwell CIPQ(SKCIPQ) if $P$ is
isomorphic to the Keedwell CIPQ, say under a mapping $\phi$.
\end{mydef}
The following results that have recently been established are of
paramount importance to prove the main result of this work.

\begin{myth}\label{1:4}(Ja\'iy\'e\d ol\'a \cite{sma14})

Let $U=(L,\oplus)$ and $V=(L,\otimes )$ be initial S-quasigroups
such that $SAUM(U)$ and $SAUM(V)$ are conjugates in $SSYM(L)$ i.e
there exists a $\psi\in SSYM(L)$ such that for any $\gamma\in
SAUM(V)$, $\gamma =\psi^{-1}\alpha\psi$ where $\alpha\in SAUM(U)$.
Then, $H_S(U)\succsim H_S(V)$ if and only if $x\delta\otimes y\gamma
=(x\beta\oplus y)\delta~\forall~x,y\in L,~\beta\in SAUM(U)$ and some
$\delta,\gamma\in SAUM(V)$.
\end{myth}

\begin{myth}\label{3:3.2}(Ja\'iy\'e\d ol\'a \cite{sma15})

The holomorph $H(L)$ of a quasigroup(loop) $L$ is a Smarandache
CIPQ(CIPL) if and only if $SAUM(L)=\{I\}$ and $L$ is a Smarandache
CIPQ(CIPL).
\end{myth}

\section{Main Results}
\begin{myth}\label{1:6}
Let $U=(L,\oplus)$ and $V=(L,\otimes )$ be initial
S-quasigroups(S-loops) that are S-isotopic under the triple of the
form $(\delta^{-1}\beta ,\gamma^{-1},\delta^{-1})$ for all $\beta\in
SAUM(U)$ and some $\delta,\gamma\in SAUM(V)$ such that their
Smarandache automorphism groups are non-trivial and are conjugates
in $SSYM(L)$ i.e there exists a $\psi\in SSYM(L)$ such that for any
$\gamma\in SAUM(V)$, $\gamma =\psi^{-1}\alpha\psi$ where $\alpha\in
SAUM(U)$. Then, $U$ is a SCIPQ(SCIPL) if and only if $V$ is a
SCIPQ(SCIPL).
\end{myth}
{\bf Proof}\\
Following Theorem~\ref{1:4}, $H_S(U)\succsim H_S(V)$. Also, by
Theorem~\ref{3:3.2}, $H_S(U)$($H_S(V)$) is a SCIPQ(SCIPL) if and
only if $SAUM(U)=\{I\}$($SAUM(V)=\{I\}$) and $U$($V$) is a
SCIPQ(SCIPL).

Let $U$ be an SCIPQ(SCIPL), then since $H_S(U)$ has a
subquasigroup(subloop) that is isomorphic to a
S-CIP-subquasigroup(subloop) of $U$ and that subquasigroup(subloop)
is isomorphic to a S-subquasigroup(subloop) of $H_S(V)$ which is
isomorphic to a S-subquasigroup(subloop) of $V$, $V$ is a
SCIPQ(SCIPL). The proof for the converse is similar.

\paragraph{Application To Cryptography}
Let the Smarandache Keedwell CIPQ be the SCIPQ $U$ in
Theorem~\ref{1:6}. Definitely, its Smarandache automorphism group is
non-trivial because as shown in Theorem~2.1 of Keedwell
\cite{phd176}, for any CIPQ, the mapping $J_\rho~:~x\to x^\rho$ is
an automorphism. This mapping will be trivial only if the
S-CIP-subquasigroup of $U$ is unipotent. For instance, in Example
2.1 of Keedwell \cite{phd176}, the CIPQ $(G,\circ )$ obtained is
unipotent because it was constructed using the cyclic group
$C_5=<c:~c^5=e>$ and defined as $a\circ b=a^3b^2$. But in
Example~2.2, the CIPQ gotten is not unipotent as a result of using
the cyclic group $C_{11}=<c:~c^{11}=e>$. Thus, the choice of a
Smarandache Keedwell CIPQ which suits our purpose in this work for a
cyclic group of order $n$ is one in which $rs=n+1$ and $r+s\ne n$.
Now that we have seen a sample for the choice of $U$, the initial
S-quasigroup $V$ can then be obtained as shown in Theorem~\ref{1:6}.
By Theorem~\ref{1:6}, $V$ is a SCIPQ.

Now, according to Theorem~\ref{1:4}, by the choice of the mappings
$\alpha ,\beta\in SAUM(U)$ and $\psi\in SSYM(L)$ to get the mappings
$\delta,\gamma$, a SCIPQ $V$ can be produced following
Theorem~\ref{1:6}. So, the secret keys for the systems are $\{\alpha
,\beta,\psi,\phi\}\equiv\{\delta,\gamma ,\phi\}$. Thus whenever a
set of information or messages is to be transmitted, the sender will
enciphere in the Smarandache Keedwell CIPQ by using specifically the
S-CIP-subquasigroup in it(as described earlier on in the
introduction) and then enciphere again with $\{\alpha
,\beta,\psi,\phi\}\equiv\{\delta,\gamma ,\phi\}$ to get a SCIPQ $V$
which is the set of encoded messages. At the receiving end, the
message $V$ is deciphered by using an inverse isotopism(i.e inverse
key of $\{\alpha ,\beta,\psi\}\equiv\{\delta,\gamma\}$) to get $U$
and then deciphere again(as described earlier on in the
introduction) to get the messages. The secret key can be changed
over time. The method described above is a double encryption and its
a double protection. It protects each piece of information(element
of the quasigroup) and protects the combined information(the
quasigroup as a whole). Its like putting on a pair of socks and
shoes or putting on under wears and clothes, the body gets better
protection. An added advantage of the use of Smarandache Keedwell
CIPQ over Keedwell CIPQ in double encryption is that the since the
S-CIP-subquasigroups of the Smarandache Keedwell CIPQ in use could
be more than one, then, the S-CIP-subquasigroups can be replaced
overtime.

\end{document}